\documentclass[a4paper,11pt]{article} 
\usepackage{amsmath,amssymb,enumerate,color}
\usepackage{amsthm,color}
\usepackage{txfonts}
\usepackage[pagewise]{lineno}
\usepackage{bigints}
\usepackage{comment}
\usepackage{graphicx}

\newcommand{\R}{\mathbb{R}}
\newcommand{\B}{\mathcal{B}}

\newcommand{\la}{\langle}
\newcommand{\ra}{\rangle}
\newcommand{\norm}[1]{\left\|#1\right\|}
\newcommand{\PV}{\text{P.V.}}
\newcommand{\abs}[1]{|#1|}
\newcommand{\ind}{\mathbf{1}}

\newtheorem{theorem}{Theorem}[section]
\newtheorem{lemma}[theorem]{Lemma}

\theoremstyle{remark}
\newtheorem{remark}{Remark}[section]
\theoremstyle{definition}

\newtheorem{definition}{Definition}[section]

\numberwithin{equation}{section}

\makeatletter
\def\@cite#1#2{[{{\bfseries #1}\if@tempswa , #2\fi}]}
\makeatother
\let\oldequation\equation
\let\oldendequation\endequation
\renewenvironment{equation}{\linenomathNonumbers\oldequation}{\oldendequation\endlinenomath}

\begin{document}
\begin{center}
\Large{{\bf
Weighted Fractional Commutator Estimate and Application to the Fractional Wave Equations with Time-dependent Damping}}
\end{center}

\vspace{5pt}

\begin{center}
Jiayun Lin
\footnote{
Department of Mathematics and Science, School of Sciences, Zhejiang Sci-Tech University, 310018, Hangzhou, China,
E-mail:\ {\tt jylin@zstu.edu.cn}},\
Masahiro Ikeda
\footnote{
Graduate School of Information Science and Technology, The University of Osaka, 1-5, Yamadaoka, Suita-shi, Osaka 565-0871, Japan/Center for Advanced Intelligence Project RIKEN, 1-4-1, Nihonbashi, Chuo-ku, Tokyo 103-0027, Japan,
E-mail:\ {\tt ikeda@ist.osaka-u.ac.jp/masahiro.ikeda@a.riken.jp}}\

\end{center}

\newenvironment{summary}{\vspace{.5\baselineskip}\begin{list}{}{%
     \setlength{\baselineskip}{0.85\baselineskip}
     \setlength{\topsep}{0pt}
     \setlength{\leftmargin}{12mm}
     \setlength{\rightmargin}{12mm}
     \setlength{\listparindent}{0mm}
     \setlength{\itemindent}{\listparindent}
     \setlength{\parsep}{0pt}
     \item\relax}}{\end{list}\vspace{.5\baselineskip}}
\begin{summary}
{\footnotesize {\bf Abstract.}
In this paper, we first establish a weighted bilinear fractional commutator estimate (Theorem~\ref{thm:weighted_commutator}) involving the spatial weight $\langle x \rangle^{2m}$. This estimate provides a powerful tool to resolve the difficulties arising from the nonlocality of the fractional Laplacian $(-\Delta)^{\sigma/2}$ and the interaction between the unbounded weights and the nonlocal operator at infinity. We then apply this result to investigate the Cauchy problem for the fractional wave equation with time-dependent damping and source nonlinearity:
\begin{equation*}
\begin{cases}
\partial_t^2 u(t,x) + (-\Delta)^{\sigma/2} u(t,x) + b(t) \partial_t u(t,x) = |u(t,x)|^p, & (t,x) \in (0,T) \times \mathbb{R}^N, \\
u(0,x) = \varepsilon u_0(x), \quad \partial_t u(0,x) = \varepsilon u_1(x), & x \in \mathbb{R}^N,
\end{cases}
\end{equation*}
where $\sigma \in (0,2)$, the initial data scale with a small parameter $\varepsilon > 0$, and the time-dependent damping coefficient $b(t)$ satisfies $b(t) \approx (1+t)^{-\beta}$ for $\beta \in [-1, 1)$. Utilizing scaling variables and constructing a weighted energy framework based on Theorem~\ref{thm:weighted_commutator}, we derive a priori differential inequalities for the total energy. For the supercritical case $p > p_c := 1 + \frac{\sigma}{N}$, we prove the small-data global existence and uniqueness of strong solutions in weighted Sobolev spaces $H^{\sigma/2,m}(\mathbb{R}^N) \times H^{0,m}(\mathbb{R}^N)$ by closing a bootstrap argument. For the subcritical ($1 < p < p_c$) and critical ($p = p_c$) cases, we establish precise lower bounds for the lifespan $T_\varepsilon$ of blow-up solutions via delicate ODE analysis. Together with the corresponding upper estimates, our results fully determine the sharp lifespan estimates for all $\beta \in [-1, 1)$, confirming that $p_c = 1 + \frac{\sigma}{N}$ is indeed the critical exponent for the fractional wave equation with time-dependent damping.
}
\end{summary}

\tableofcontents

\section{Introduction}

During the past few decades, the fractional Laplacian \((-\Delta)^{\sigma/2}\) with \(\sigma \in (0, 2)\) has become a fundamental tool for modeling anomalous diffusion and nonlocal phenomena. It frequently appears in various partial differential equations, such as fractional diffusion-reaction equations \cite{YA}, quasi-geostrophic models \cite{CW}, and fractional Schrödinger equations \cite{Las}. Mathematically, it can be defined through the Fourier transform by
\[
\mathcal{F}\left[(-\Delta)^{\sigma/2} f\right](\xi) = |\xi|^\sigma \hat{f}(\xi),
\]
or equivalently as the singular integral
\[
(-\Delta)^{\sigma/2} f(x) = C_{N,\sigma} \text{P.V.} \int_{\mathbb{R}^N} \frac{f(x) - f(y)}{|x-y|^{N+\sigma}} \, dy,
\]
where $C_{N,\sigma} > 0$ is a normalization constant. 

From a stochastic viewpoint, for $0<\sigma<2$, the operator $(-\Delta)^{\sigma/2}$ is the negative of the infinitesimal generator of a symmetric $\sigma$-stable L\'evy process \cite{Lis}. This process is characterized by heavy-tailed jump distributions, which give rise to increments with infinite variance---a defining feature of L\'evy flights. As a result, particles perform long-range jumps, and the motion exhibits superdiffusive scaling, with the typical displacement growing as $t^{1/\sigma}$. 

In the study of nonlinear PDEs involving the fractional Laplacian, commutator estimates serve as a key analytic tool. In particular, they play a decisive role in the fixed point arguments used to establish local well-posedness in Sobolev spaces for nonlinear dispersive PDEs. The classical Kato--Ponce commutator estimate \cite{KP} states that for \(1<p<\infty\) and \(\sigma>0\), and for \(f,g \in \mathcal{S}(\mathbb{R}^N)\),
\[
\| J^\sigma(fg) - f J^\sigma g \|_{L^p(\mathbb{R}^N)} 
\lesssim \| J^\sigma f \|_{L^p(\mathbb{R}^N)} \| g \|_{L^\infty(\mathbb{R}^N)} 
+ \| \nabla f \|_{L^\infty(\mathbb{R}^N)} \| J^{\sigma-1} g \|_{L^p(\mathbb{R}^N)},
\]
where \(J^\sigma := (1-\Delta)^{\sigma/2}\) and \(\nabla = (\partial_1,\dots,\partial_N)\). This classical estimate and its variants have been extended to various geometric settings and weighted spaces (see \cite{FR} for example). However, in the present work we focus on a particular weighted bilinear commutator estimate involving the polynomial weight \(|x|^{2m}\), which arises naturally in the study of the fractional wave equations with time-dependent damping. Our first result is the following.

\begin{theorem}
\label{thm:weighted_commutator}
Let \(0 < \sigma < 2\) and \(0 < \varepsilon < \sigma\). Consider the weight function \(w(x) = |x|^{2m}\). Assume that 
\[
m > \frac{N}{2} + 1, \qquad \text{and} \qquad m + \varepsilon < \frac{N}{2} + \sigma + 1.
\]
If \(u, U \in L^2(\R^N)\) satisfy the zero-mean conditions
\[
\int_{\R^N} u(x)\,dx = 0 \qquad \text{and} \qquad \int_{\R^N} U(x)\,dx = 0,
\]
then there exists a constant \(C > 0\), depending only on \(N, \sigma, m,\) and \(\varepsilon\), such that the bilinear commutator integral \(\mathcal{B}(u,U) := \int_{\R^N} [\Lambda^\sigma, w]u(x)U(x)\,dx\) satisfies the following estimate:
\begin{equation}\label{est:fraccommu}
\begin{aligned}
&\left| \int_{\R^N} [\Lambda^\sigma, w]u(x)U(x)\,dx \right|\\
&\qquad \le
C \norm{\la x\ra^{m-\varepsilon}u}_{L^2} \norm{\la x\ra^m U}_{L^2}
+ C \norm{\la x\ra^m u}_{L^2} \norm{\la x\ra^{m-\varepsilon}U}_{L^2}.
\end{aligned}
\end{equation}
Here, \(\Lambda^\sigma = (-\Delta)^{\sigma/2}\) denote the fractional Laplacian of order \(\sigma\).
\end{theorem}

\begin{remark}
The weighted bilinear fractional commutator estimate \eqref{est:fraccommu} naturally extends to the regularized Japanese bracket weight $w(x) = \langle x \rangle^{2m}$ with the same argument.
\end{remark}

We now apply Theorem~\ref{thm:weighted_commutator} to the Cauchy problem of the fractional wave equations with time-dependent damping:
\begin{equation}\label{main}
\begin{cases}
\partial_t^2 u(x,t) + (-\Delta)^{\sigma/2} u(x,t) + b(t) \partial_t u(x,t) = |u(x,t)|^p, & (x,t) \in \mathbb{R}^N \times (0,T), \\
u(x,0) = \varepsilon u_0(x), \quad \partial_t u(x,0) = \varepsilon u_1(x), & x \in \mathbb{R}^N,
\end{cases}
\end{equation}
where $T\in (0,\infty]$ is the existence time of the function $u$, $u:[0,T)\times \R^N$ is an unknown function. Here $b:[0,\infty)\rightarrow \mathbb{R}_{>0}$ is a given positive function, satisfying
\begin{equation}\label{b}
    b_1(1+t)^{-\beta}\leq b(t)\leq b_2(1+t)^{-\beta},\quad \quad \mbox{and}\quad \quad |b'(t)|\leq b_3(1+t)^{-1}b(t),
\end{equation}
for $\beta\in \R$ and some constants $b_1, b_2, b_3>0$ and $b(0)=1$. 

Our aim is to establish global existence of solutions in the supercritical case and sharp lifespan estimates in the subcritical and critical cases. The weighted commutator estimate proved above plays a central role. Indeed, when deriving weighted \(L^2\)-type estimates for the solution, one is naturally led to commute the fractional Laplacian \(\Lambda^\sigma = (-\Delta)^{\sigma/2}\) with the polynomial weight \(|x|^{2m}\). Theorem~\ref{thm:weighted_commutator} provides the required bilinear control of the resulting commutator, which is essential for closing the bootstrap argument in the supercritical case and for obtaining the lower lifespan bounds in the remaining cases.

The long-time behavior of solutions to the Cauchy problem for semilinear damped wave equations has been extensively studied. In the classical case $\sigma=2$ and $b(t)=1$, Todorova and Yordanov \cite{TY1, TY2}, and Zhang \cite{Z} established the critical exponent $p_c = 1 + \frac{2}{N}$, proving global existence for small data in the supercritical case ($p_c <\rho< \frac {N+2}{[N-2]_+}$) and blow-up in the subcritical and critical cases ($1< p \le p_c$). The sharp lifespan $T(\varepsilon)$ for blow-up solutions is bounded by $C \varepsilon^{-(\frac{1}{p-1}-\frac{N}{2})^{-1}}$ for $1 < p < p_c$, and by $\exp(C\varepsilon^{-(p-1)})$ for $p = p_c$, see \cite{LZ, N2003, IO, LaiZhou}.

For time-dependent damping $b(t)=(1+t)^{-\beta}$ with $-1<\beta<1$, the critical exponent remains $p_c=1+\frac{2}{N}$, see \cite{NLZ, N}. The sharp lifespan estimate in the subcritical case is given by
\[
T(\varepsilon)\approx C\varepsilon^{-\left[\left(\frac{1}{p-1}-\frac{N}{2}\right)(1+\beta)\right]^{-1}},
\]
see \cite{FIW2019, IW}. In the critical case, the lifespan remains exponentially bounded:
\[
T(\varepsilon)\approx \exp\bigl(C\varepsilon^{-(p-1)}\bigr),
\]
as shown in \cite{IO,ISW, IW}. When $\beta=-1$, the lifespan exhibits exponential and double-exponential bounds for the subcritical and critical cases:
\[
T(\varepsilon)\approx
\begin{cases}
C\exp\Bigl(\varepsilon^{-\left(\frac{1}{p-1}-\frac{N}{2}\right)^{-1}}\Bigr), & 1<p<p_c,\\[0.5ex]
\exp\bigl(\exp(C\varepsilon^{-(p-1)})\bigr), & p=p_c,
\end{cases}
\]
see \cite{FIW2019,II}. For more general effective damping, we refer to \cite{DLR,DL,ISW}.

For the fractional damped equation with constant damping, the critical exponent is $p_c=1+\frac{\sigma}{N}$. Dao and Reissig \cite{DR} proved blow-up for $1<p<p_c$ when $\sigma\geq2$, while Fujiwara--Ikeda--Wakasugi \cite{FIW2020} obtained small-data global existence for $p>p_c$ in suitable weighted Sobolev spaces. For weakly coupled systems, see \cite{DP}. Recently, for time-dependent damping and $\sigma>2$, Aslan and Dao \cite{AD} proved global well-posedness in the supercritical case and derived sharp lifespan estimates in the subcritical and critical cases.

Before stating the main results, we first recall the weighted Sobolev spaces used in this paper. For $s\in \R$ and $m\ge 0$, we denote the Sobolev space $H^{s,m}(\R^N)$ with the norm
\begin{equation*}
\|f\|_{H^{s,m}}:= \left(\int_{\R^N}\la x \ra^{2m}(\langle \nabla \rangle^s f)^2dx\right)^{1/2},
\end{equation*}
where $\langle \cdot \rangle :=(1+\vert \cdot \vert ^2)^{1/2}$.

With this notation, we are now ready to define the strong solution and its lifespan for the Cauchy problem \eqref{main}.
\begin{definition}[Strong solution and its lifespan]
Let $N\in \mathbb{N}$, $\sigma\in (0,2)$, $(u_0,u_1)\in H^{\sigma,m}(\R^N)\times H^{\sigma/2,m}(\R^N)$ with $m\ge 0$ and $T>0$. We say that $u$ is a strong solution to the Cauchy problem \eqref{main} on $[0,T)$ if $u\in C^2([0,T);H^{0,m}(\R^N))\cap C^1([0,T);H^{\sigma/2,m}(\R^N))\cap C([0,T);H^{\sigma,m}(\R^N))$ satisfies the initial condition $u(0)=\varepsilon u_0$ and $\partial_tu(0)= \varepsilon u_1$ and satisfies the equation
\[
\partial_t^2u+(-\Delta)^{\sigma/2}u+b(t)\partial_tu=f(u)
\]
in the sense of $C([0,T);H^{0,m}(\R^N))$. Moreover, the lifespan \(T(\varepsilon)\) of a solution \(u\) for the Cauchy problem (\ref{main}) is defined by
\[
T(\varepsilon) = \sup\{T > 0 \mid u \text{ is a strong solution for \eqref{main} on } [0, T]\}.
\]
\end{definition}

By the local existence theory of strong solutions established in \cite{LI}, for any initial data $(u_0, u_1) \in H^{\sigma,m}(\mathbb{R}^N) \times H^{\sigma/2,m}(\mathbb{R}^N)$, there exists a maximal existence time $T(\varepsilon) > 0$ such that the Cauchy problem \eqref{main} admits a unique strong solution on $[0, T(\varepsilon))$. Based on this local existence, by establishing uniform global energy estimates to prevent blow-up, we obtain the global existence of solutions in the supercritical case $p > 1 + \frac{\sigma}{N}$, as stated in the following theorem.

\begin{theorem}[Global existence]
\label{globalexistence}
Let $\sigma \in (0,2)$, $\beta \in [-1,1)$, $1+\frac{\sigma}{N} < p <\infty$ for $0 < N \le \sigma$ and $1+\frac{\sigma}{N} < p < \frac N{N-\sigma}$ for $N>\sigma$. Assume that the initial data $(u_0,u_1)\in H^{\sigma/2,m}\times H^{0,m}$, where $m=1\ (N=1)$ and $1+\frac{N}{2}<m<1+\frac{N}{2}+\sigma\ (N \ge 2)$. There exists $\varepsilon_0 > 0$ such that for every $\varepsilon \in (0,\varepsilon_0]$, there exists a unique solution 
\[
u \in C([0,\infty);H^{\sigma/2,m}(\R^N))\cap C^1([0,\infty);H^{0,m}(\R^N))
\]
for the Cauchy problem \eqref{main}.
\end{theorem}
Let us introduce the main result on the lower bound of the lifespan in the subcritical and critical cases.

\begin{theorem}[Lower bound of the lifespan in the subcritical case]
\label{blowupsub}
Let $\sigma \in (0,2)$, $\beta \in [-1,1)$, and $1<p<1+\frac{\sigma}{N} $. Assume that the initial data $(u_0,u_1)\in H^{\sigma/2,m}\times H^{0,m}$, where $m=1\ (N=1)$ and $1+\frac{N}{2}<m<1+\frac{N}{2}+\sigma\ (N \ge 2)$. Then, there exists a positive constant $\varepsilon_1>0$, such that for every $\varepsilon \in (0,\varepsilon_1]$, the lower bound of the lifespan of the strong solution is estimated by
\begin{equation}
T\ge
\begin{cases}
C \varepsilon^{-\frac 1{\left(\frac 1{p-1}-\frac N\sigma\right)(1+\beta)}}, &\quad \mbox{for}\quad -1<\beta\leq 1,\\
\exp{\left(C\varepsilon^{-\frac{1}{\frac 1{p-1}-\frac N\sigma}}\right)}, &\quad \mbox{for}\quad \beta= -1,
\end{cases}
\end{equation}
for some constant $C$.
\end{theorem}

And in the critical case, we have
\begin{theorem}[Lower bound of the lifespan in the critical case]
\label{blowupcritical}
Let $\sigma \in (0,2)$, $\beta \in [-1,1)$, and $p= 1+\frac \sigma N$. Assume that the initial data $(u_0,u_1)\in H^{\sigma/2,m}\times H^{0,m}$, where $m=1\ (N=1)$ and $1+\frac{N}{2}<m<1+\frac{N}{2}+\sigma\ (N \ge 2)$. Then, there exists a positive constant $\varepsilon_2>0$, such that for every $\varepsilon \in (0,\varepsilon_2]$, the lower bound of the lifespan of the strong solution is estimated by
\begin{equation}
    T\ge
    \begin{cases}
        \exp{(C\varepsilon^{-(p-1)})},\quad &\mbox{if}\quad \beta\in (-1,1),\\
        \exp{(\exp{(C\varepsilon^{-(p-1)}}))},\quad &\mbox{if}\quad \beta=-1,
    \end{cases}
\end{equation}
for some constant $C$.
\end{theorem}

\begin{remark}
Together with the upper bound derived in \cite{LI}, our lower bound establishes the sharp estimate for the lifespan.
\begin{itemize}
    \item[\rm(i)] \textbf{Subcritical case:} If $1 < p < 1 + \frac{\sigma}{N}$, 
    \begin{equation}
    T\approx
    \begin{cases}
        C \varepsilon^{-\frac 1{\left(\frac 1{p-1}-\frac N\sigma\right)(1+\beta)}}, &\quad \mbox{for}\quad -1<\beta\leq 1,\\
        \exp{\left(C\varepsilon^{-\frac{1}{\frac 1{p-1}-\frac N\sigma}}\right)}, &\quad \mbox{for}\quad \beta= -1,
    \end{cases}
\end{equation}
    \item[\rm(ii)] \textbf{Critical case:} If $p = 1 + \frac{\sigma}{N}$, 
    \begin{equation}
        T\approx
        \begin{cases}
            \exp{(C\varepsilon^{-(p-1)}},\quad &\mbox{if}\quad \beta\in (-1,1),\\
            \exp{(\exp{(C\varepsilon^{-(p-1)}})},\quad &\mbox{if}\quad \beta=-1,
        \end{cases}
    \end{equation}
\end{itemize}
\end{remark}

\begin{remark}
As a consequence of Theorem~\ref{globalexistence}, Theorems~\ref{blowupsub}--\ref{blowupcritical}, and the upper bounds derived in \cite{LI}, the critical exponent for the Cauchy problem \eqref{main} is precisely 
\[
p_c = 1 + \frac{\sigma}{N}.
\]
Namely, the problem admits global strong solutions for $p > p_c$ with small initial data, whereas finite-time blow-up occurs for $1 < p \le p_c$ with the sharp estimate of the lifespan $T(\varepsilon)$.
\end{remark}

The main proof strategy of Theorem \ref{globalexistence}-Theorem \ref{blowupcritical} relies on a weighted fractional commutator estimate. The primary difficulty in establishing this estimate arises from the mismatch between the non-local nature of the fractional Laplacian $(-\Delta)^{\sigma/2}$ and the unboundedness of the polynomial weights. Unlike integer-order derivatives, the fractional operator lacks a standard Leibniz rule, meaning that the commutator cannot be simply reduced to lower-order local terms. Furthermore, the singular integral representation of the fractional operator strongly interacts with the growth of the weights at infinity, creating difficulty in obtaining suitable weighted norm estimates. To overcome these analytical difficulties, our proof relies on a delicate region decomposition strategy combined with precise kernel estimates. We carefully partition the domain of integration in order to separate the local singularity of the fractional kernel from the long-range interactions generated by the polynomial weights. In addition, we adapt and extend classical singular integral results and properties of Riesz potentials to our weighted framework, which allows us to control the remainder terms produced by the nonlocal interactions. Under this framework, we adopt the scaling variables introduced by Wakasugi \cite{Wakasugi2017} and employ a weighted energy method to conduct a rigorous proof. Briefly, the approach first transforms the original equation via scaling variables and decomposes the solution in such a way that the remainder has zero mean. Subsequently, by constructing weighted energies and incorporating our commutator estimates, we derive a crucial ordinary differential inequality for the total energy. This fundamental inequality enables us to close the global energy estimates in the supercritical case via a bootstrap argument, and to explicitly compute the sharp lower bounds of the lifespan in the subcritical and critical cases.  

The rest of this paper is organized as follows. In Section 2, we establish the weighted fractional commutator estimate, along with weighted bounds for the Bessel potential operator. Section 3 is devoted to the construction of the weighted energy framework, where we introduce scaling variables, decompose the solutions, and derive the key a priori differential inequalities for the constructed energies. In Section 4, we prove the global existence of solutions in the supercritical case utilizing the bootstrap argument. Finally, in Section 5, we derive the lower bounds of the lifespan in the subcritical and critical cases by analyzing the integral inequality. 

\section{Fractional commutator estimates}

In this section we prove Theorem~\ref{thm:weighted_commutator}, which is the key analytical ingredient for the energy estimates in later sections.  
The main difficulty is that the fractional Laplacian is nonlocal and does not obey a standard Leibniz rule, so the commutator with the unbounded weight \(w(x)=|x|^{2m}\) cannot be handled by purely local arguments. To overcome this, we decompose the double integral defining the bilinear form into several regions by means of a smooth partition of unity. Roughly speaking, the decomposition distinguishes between the local region where \(|x-y|\lesssim 1\), the region where \(|x|\sim |y|\), and the far-field regions where one variable dominates the other. In the local and comparable regions, the weight difference \(w(x)-w(y)\) can be controlled directly by kernel estimates and Young's inequality. In the far-field regions, the zero-mean conditions on \(u\) and \(U\) play a crucial role. They allow us to subtract a suitable reference value and gain the additional decay needed to make the corresponding integrals finite. Combining the estimates from all regions yields the desired weighted bilinear bound.

We now turn to the detailed proof.

\begin{proof}[Proof of Theorem~\ref{thm:weighted_commutator}]
For the given function \(w\), the commutator \([\Lambda^s, w]\) acting on a function \(u\) is defined by
\begin{equation}
[\Lambda^s, w]u = \Lambda^s(wu) - w\Lambda^s u.
\end{equation}
Note that 
\[
\Lambda^s f(x) = c_{N,s} \,\PV \int_{\R^N} \frac{f(x)-f(y)}{|x-y|^{N+s}}\,dy,
\]
where \(c_{N,s}\) is a normalization constant depending only on the dimension \(N\) and \(s\). Thus, the commutator could be expressed as
\begin{equation}
\begin{aligned}
\relax [\Lambda^s, w]u & = \Lambda^s(wu) - w\Lambda^s u\\
& = c_{N,s} \,\PV \int_{\R^N} \frac{w(x)u(x)-w(y)u(y)}{|x-y|^{N+s}}\,dy\\
&\qquad - w(x)\left( c_{N,s} \,\PV \int_{\R^N} \frac{u(x)-u(y)}{|x-y|^{N+s}}\,dy \right) \\
& = c_{N,s} \,\PV \int_{\R^N} \frac{(w(x)-w(y))u(y)}{|x-y|^{N+s}}\,dy.
\end{aligned}
\end{equation}
Therefore, introducing the kernel
\[
K(x,y) := \frac{w(x)-w(y)}{|x-y|^{N+s}},
\]
the bilinear form can be written as
\begin{equation}
\B(u, U) = c_{N,s} \int_{\R^N} \int_{\R^N} K(x,y)\,u(y)\,U(x)\,dy\,dx.
\end{equation}

We split the integration domain into four exactly separated parts by a smooth partition of unity, and estimate each part separately.

Choose smooth cut-off functions $\eta_0,\eta_1\in C^\infty([0,\infty))$ such that
\[
0\le \eta_0,\eta_1\le 1,\quad
\eta_0(t)+\eta_1(t)=1,
\]
\[
\eta_0(t)=1\;(0\le t\le 1),\qquad
\eta_0(t)=0\;(t\ge 2).
\]
Thus, $\eta_1$ is supported where $t\ge 1$.

Next, choose $\chi_X,\chi_Y,\chi_C\in C^\infty((0,\infty))$ satisfying
\[
\chi_X(r)+\chi_Y(r)+\chi_C(r)=1\quad (r>0),
\]
\[
\chi_X(r)=1\;(0<r\le 1/3),\qquad
\chi_X(r)=0\;(r\ge 1/2),
\]
\[
\chi_Y(r)=1\;(r\ge 3),\qquad
\chi_Y(r)=0\;(r\le 2).
\]
Hence, $\chi_C$ is supported in $[1/3,3]$.

For $x,y\in\R^N$, define
\[
r(x,y):=\frac{\la y\ra}{\la x\ra}.
\]
Noting that $\eta_0+\eta_1=1$ and $\chi_X+\chi_C+\chi_Y=1$, we obtain the exact identity
\[
1 = \eta_0(|x-y|) + \eta_1(|x-y|)\Bigl[\chi_X(r(x,y)) + \chi_C(r(x,y)) + \chi_Y(r(x,y))\Bigr].
\]
Inserting this partition of unity into the double integral for $\B(u,U)$ yields
\begin{equation}
\B = \B_{\rm loc} + \B_X + \B_C + \B_Y,
\end{equation}
where
\begin{align*}
\B_{\rm loc} &:= c_{N,s}\int_{\R^N} \int_{\R^N} \eta_0(|x-y|)\,K(x,y)\,u(y)\,U(x)\,dy\,dx,\\
\B_X        &:= c_{N,s}\int_{\R^N} \int_{\R^N} \eta_1(|x-y|)\,\chi_X(r(x,y))\,K(x,y)\,u(y)\,U(x)\,dy\,dx,\\
\B_C        &:= c_{N,s}\int_{\R^N} \int_{\R^N} \eta_1(|x-y|)\,\chi_C(r(x,y))\,K(x,y)\,u(y)\,U(x)\,dy\,dx,\\
\B_Y        &:= c_{N,s}\int_{\R^N} \int_{\R^N} \eta_1(|x-y|)\,\chi_Y(r(x,y))\,K(x,y)\,u(y)\,U(x)\,dy\,dx.
\end{align*}
Then, we estimate each term separately.

{\bf The term $\B_{\rm loc}$.}
On the support of $\eta_0(|x-y|)$ we have $|x-y|\le 2$, and consequently $\la x\ra\sim\la y\ra$.
By the mean value theorem and the bounds on $w$,
\[
\abs{w(x)-w(y)}
\le C\,|x-y|\bigl(\la x\ra^{2m-1}+\la y\ra^{2m-1}\bigr).
\]
Since $\la x\ra\sim\la y\ra$ and $0<\varepsilon<1$,
\[
\la x\ra^{2m-1}\le C\la x\ra^m\la y\ra^{m-\varepsilon},\quad \text{and}\quad
\la y\ra^{2m-1}\le C\la x\ra^m\la y\ra^{m-\varepsilon}.
\]
Thus,
\begin{equation}
\abs{K(x,y)}\,\eta_0(|x-y|)
\le C\,
\frac{\ind_{\{|x-y|\le 2\}}}{|x-y|^{N+s-1}}
\la x\ra^m\la y\ra^{m-\varepsilon}.
\end{equation}

Since $s<1$, 
\[
\int_{\R^N} \ind_{\{|z|\le 2\}}\,|z|^{-N-s+1}\,dz = C\int_0^2 r^{-s}\,dr <\infty.
\]
Applying Young's inequality, we obtain
\begin{equation}
\begin{aligned}
\abs{\B_{\rm loc}} &\le C\int_{\R^N} \int_{\R^N} \frac{\ind_{\{|x-y|\le 2\}}}{|x-y|^{N+s-1}} \la x\ra^m\la y\ra^{m-\varepsilon} \abs{u(y)}\,\abs{U(x)}\,dy\,dx\\
&\le C\, \norm{\la x\ra^{m-\varepsilon}u}_{L^2}\norm{\la x\ra^mU}_{L^2}.
\end{aligned}
\end{equation}
This gives the desired estimate for $\B_{\rm loc}$.

{\bf The term $\B_C$.}
On the support of $\chi_C(r(x,y))$, we have $\la x\ra\sim\la y\ra$, and $\eta_1(|x-y|)$ implies $|x-y|\ge 1$.
In this region, we claim the bound
\begin{equation}
\abs{w(x)-w(y)}
\le C \la x\ra^m\la y\ra^{m-\varepsilon} \,|x-y|^\varepsilon.
\end{equation}

To prove the claim, consider two subregions.
\begin{itemize}
\item If $|x-y|\ge c\la y\ra$, then $\la y\ra^\varepsilon\le C|x-y|^\varepsilon$.  Since
$$\abs{w(x)-w(y)}\le C(\la x\ra^{2m}+\la y\ra^{2m})\le C\la x\ra^m\la y\ra^m,$$
we obtain
\[
\abs{w(x)-w(y)}\le C\la x\ra^m\la y\ra^{m-\varepsilon}|x-y|^\varepsilon.
\]
\item If $1\le |x-y|<c\la y\ra$, the mean value theorem gives
\[
\abs{w(x)-w(y)}\le C|x-y|\,\la x\ra^{2m-1}\quad(\text{since }\la x\ra\sim\la y\ra).
\]
Since $\varepsilon-1<0$ and $1\le|x-y| < c\la y\ra \le C\la x\ra$, we have 
\[
\la x\ra^{\varepsilon-1}\le C|x-y|^{\varepsilon-1},
\]
which implies $\la x\ra^{2m-1}\le C\la x\ra^m\la y\ra^{m-\varepsilon}\la x\ra^{\varepsilon-1}$. Hence
\[
\abs{w(x)-w(y)}\le C\la x\ra^m\la y\ra^{m-\varepsilon}|x-y|^\varepsilon.
\]
\end{itemize}
Combining the estimates obtained in two subregions, we could derive the claim.

Consequently, on the support of $\eta_1\chi_C$,
\begin{equation}
\abs{K(x,y)}\le C\,
\frac{\la x\ra^m\la y\ra^{m-\varepsilon} }{|x-y|^{N+s-\varepsilon}}.
\end{equation}
Since $0<\varepsilon<s$, we have $-N-s+\varepsilon<-N$, hence $\ind_{\{|z|\ge 1\}}\,|z|^{-N-s+\varepsilon} \in L^1(\R^N)$.  Repeating the Young inequality argument as for the first term, we obtain
\begin{equation}
\abs{\B_C}
\le C \norm{\la x\ra^{m-\varepsilon}u}_{L^2}\norm{\la x\ra^mU}_{L^2}.
\end{equation}

{\bf The term $\B_X$.}
Here $\chi_X(r(x,y))$ is active, which enforces $\la y\ra\le\frac12\la x\ra$.
For fixed $x$, set
\[
\Phi_x(y) := \eta_1(|x-y|)\,\chi_X\!\Bigl(\frac{\la y\ra}{\la x\ra}\Bigr)\,K(x,y).
\]
Then
\[
\B_X = c_{N,s}\int_{\R^N} U(x)\Bigl(\int_{\R^N} \Phi_x(y)\,u(y)\,dy\Bigr)dx.
\]

Because 
\[
\int_{\R^N} u(x)\,dx = 0,
\]
we may write
\begin{equation}
\int_{\R^N}\Phi_x(y)\,u(y)\,dy = \int_{\R^N}\bigl(\Phi_x(y)-\Phi_x(0)\bigr)u(y)\,dy.
\end{equation}
Notice that when $y=0$, $\la y\ra=1$, and $\chi_X(1/\la x\ra)>0$ implies $\la x\ra\ge2$, hence $|x|\ge\sqrt3>1$.  Thus, $|x|>1$, which implies that 
\[
\Phi_x(0) = \eta_1(|x|)\,\chi_X\!\Bigl(\frac{1}{\la x\ra}\Bigr)\,\frac{w(x)-w(0)}{|x|^{N+s}}
\] 
is completely free of singularities.

By the mean value theorem,
\begin{equation}
\abs{\Phi_x(y)-\Phi_x(0)}\le |y|\sup_{0\le\theta\le1}\abs{\nabla_y\Phi_x(\theta y)}.
\end{equation}

We now estimate $\nabla_y\Phi_x$. Direct calculation gives
\[
\nabla_y\Phi_x = (\nabla_y\eta_1)\,\chi_X K \;+\; \eta_1(\nabla_y\chi_X)K \;+\; \eta_1\chi_X\,\nabla_yK.
\]
We examine each term separately.

For the first term, notice that $\nabla_y\eta_1(|x-y|)$ is non-zero only when $1\le|x-y|\le2$ by the definition of $\eta$. With the restriction $\la y\ra\le\frac12\la x\ra$ from $\chi_X$, we obtain
\[
|x| \le |y| + |x-y| \le \la y\ra + 2 \le \frac12\la x\ra + 2,
\]
which yields $\la x\ra \le 6$. Besides, the kernel $K(x,y)$ satisfies
\[
\abs{K(x,y)} = \frac{\abs{w(x)-w(y)}}{|x-y|^{N+s}} \le C\la x \ra^{2m-N-s},
\]
since $\la y\ra \le \frac12\la x\ra$ and $|x-y|\sim\la x\ra$ on the support of $\Phi_x(y)$. Combining these two estimates, we get
\[
\abs{(\nabla_y\eta_1)\,\chi_X K}\le C \la x \ra^{2m-1-N-s},
\]
since $\la x \ra^{-1} \ge \frac 16$.

For the second term, since $\nabla_y\chi_X(\la y\ra/\la x\ra) = \chi_X'(\cdot)\,\frac{y}{\la y\ra}\la x\ra^{-1}$, its size is at most $C\la x\ra^{-1}$.  Hence
\[
\abs{\eta_1(\nabla_y\chi_X)K}\le C\la x\ra^{-1}\la x\ra^{2m-N-s}=C\la x\ra^{2m-1-N-s}.
\]
Then, for the last term, differentiating $K(x,y)$ w.r.t. $y$, we have
\[
\nabla_yK(x,y)= -\frac{\nabla w(y)}{|x-y|^{N+s}} + (w(x)-w(y))\nabla_y\bigl(|x-y|^{-N-s}\bigr).
\]
Using the bounds on $w$ and $\nabla w$, together with the properties that $|x-y|\sim\la x\ra$ and $\la y\ra \le \frac12\la x\ra$ on the support of $\Phi_x(y)$, we obtain
\[
\abs{\nabla_yK(x,y)}
\le C\frac{\la y\ra^{2m-1}}{\la x\ra^{N+s}} + C\frac{\la x\ra^{2m}}{\la x\ra^{N+s+1}}
\le C\la x\ra^{2m-1-N-s}.
\]
Since the cut-off functions are bounded, we get
\[
\abs{\eta_1\chi_X\,\nabla_yK}\le C \la x \ra^{2m-1-N-s},
\]
Hence, globally on the support,
\begin{equation}
\abs{\nabla_y\Phi_x(y)}\le C\la x\ra^{2m-1-N-s}.
\end{equation}
Therefore,
\begin{equation}
\abs{\int_{\R^N}\Phi_x(y)u(y)\,dy}
\le C\la x\ra^{2m-1-N-s}\int_{\R^N}|y|\,|u(y)|\,dy.
\end{equation}
By the Cauchy--Schwarz inequality,
\[
\int_{\R^N}|y|\,|u(y)|\,dy
\le \Bigl(\int_{\R^N}|y|^2\la y\ra^{-2m}\,dy\Bigr)^{1/2}
\norm{\la y\ra^m u}_{L^2}.
\]
The weight integral is finite precisely because $m>\frac N2+1$.  Thus,
\begin{equation}
\abs{\int_{\R^N}\Phi_x(y)u(y)\,dy}
\le C\la x\ra^{2m-1-N-s}\norm{\la y\ra^m u}_{L^2}.
\end{equation}
Substituting this into $\B_X$ yields
\begin{equation}
\abs{\B_X}
\le C\norm{\la x\ra^m u}_{L^2}
\int_{\R^N} \la x\ra^{2m-1-N-s}\,\abs{U(x)}\,dx.
\end{equation}
Applying the Cauchy--Schwarz inequality with weight $\la x\ra^{m-\varepsilon}$,
\begin{align*}
\int_{\R^N} \la x\ra^{2m-1-N-s}\abs{U(x)}\,dx
&\le \Bigl(\int_{\R^N} \la x\ra^{2(m+\varepsilon-1-N-s)}\,dx\Bigr)^{1/2}
\norm{\la x\ra^{m-\varepsilon}U}_{L^2}.
\end{align*}
The integral converges exactly when $2(m+\varepsilon-1-N-s) < -N$, i.e. $m+\varepsilon < \frac N2 + s + 1$, which holds by the assumption.  Consequently,
\begin{equation}
\abs{\B_X}
\le C\,\norm{\la x\ra^m u}_{L^2}\,
\norm{\la x\ra^{m-\varepsilon}U}_{L^2}.
\end{equation}

{\bf The term $\B_Y$.}
This case is completely symmetric to the previous one.  For fixed $y$, define
\begin{equation}
\Psi_y(x) := \eta_1(|x-y|)\,\chi_Y\!\Bigl(\frac{\la y\ra}{\la x\ra}\Bigr)\,K(x,y).
\end{equation}
Because 
\[
\int_{\R^N} U(x)\,dx = 0,
\]
we write
\begin{equation}
\int_{\R^N}\Psi_y(x)U(x)\,dx
= \int_{\R^N}\bigl(\Psi_y(x)-\Psi_y(0)\bigr)U(x)\,dx.
\end{equation}
On the support, we have $\la x\ra\le \frac 12 \la y\ra$ and $|x-y|\sim\la y\ra$.  Repeating the derivative estimates yields
\begin{equation}
\abs{\nabla_x\Psi_y(x)}\le C\la y\ra^{2m-1-N-s}.
\end{equation}
The same reasoning as before gives
\begin{equation}
\abs{\int_{\R^N}\Psi_y(x)U(x)\,dx}
\le C\la y\ra^{2m-1-N-s}\norm{\la x\ra^m U}_{L^2}.
\end{equation}
Inserting this into the expression for $\B_Y$,
\begin{equation}
\abs{\B_Y}
\le C\norm{\la x\ra^m U}_{L^2}
\int_{\R^N} \la y\ra^{2m-1-N-s}\,\abs{u(y)}\,dy.
\end{equation}
By the Cauchy--Schwarz inequality with weight $\la y\ra^{m-\varepsilon}$ and the same condition $m+\varepsilon<\frac N2+s+1$, we conclude
\begin{equation}
\abs{\B_Y}
\le C\,\norm{\la x\ra^m U}_{L^2}\,
\norm{\la y\ra^{m-\varepsilon}u}_{L^2}
= C\,\norm{\la x\ra^{m-\varepsilon}u}_{L^2}
\norm{\la x\ra^m U}_{L^2}.
\end{equation}

Summing the estimates for $\B_{\rm loc}$, $\B_C$, $\B_X$, and $\B_Y$, we obtain
\begin{equation}
\abs{\B(u,U)}
\le C\Bigl(
\norm{\la x\ra^{m-\varepsilon}u}_{L^2}\norm{\la x\ra^mU}_{L^2}
+
\norm{\la x\ra^mu}_{L^2}\norm{\la x\ra^{m-\varepsilon}U}_{L^2}
\Bigr),
\end{equation}
which completes the proof.
\end{proof}

\section{The weighted energy framework}

Following the strategy of Wakasugi \cite{Wakasugi2017} for the classical damped wave equation, we adapt the scaling variables technique to the fractional Laplacian case. The main difficulty here is the treatment of the nonlocal operator $(-\Delta)^{\sigma/2}$, which requires new commutator identities and a more involved energy construction.

We introduce the change of variables
\begin{equation}\label{transformation}
s = \log(B(t)+1),\qquad y = (B(t)+1)^{-\frac1\sigma}x,
\end{equation}
where
\[
B(t)=\int_0^t \frac{1}{b(\tau)}\,d\tau .
\]
Then, define new unknown functions $v,w$ by
\begin{equation}\label{uvw}
\begin{aligned}
u(x,t) &= (B(t)+1)^{-\frac N\sigma}\,v\bigl((B(t)+1)^{-\frac1\sigma}x,\;\log(B(t)+1)\bigr),\\
u_t(x,t) &= b(t)^{-1}(B(t)+1)^{-\frac N\sigma-1}\,w\bigl((B(t)+1)^{-\frac1\sigma}x,\;\log(B(t)+1)\bigr).
\end{aligned}
\end{equation}
The direct calculation shows that the pair $(v,w)$ solves the system
\begin{equation}\label{system-vw}
\begin{cases}
v_s - \displaystyle\frac{1}{\sigma}\,y\cdot\nabla_y v - \frac N\sigma\,v = w,\\[10pt]
\displaystyle\frac{e^{-s}}{b(t)^2}\Bigl(w_s - \frac{y}{\sigma}\cdot\nabla_y w - \Bigl(\frac N\sigma+1\Bigr)w\Bigr) + w = -(-\Delta_y)^{\sigma/2}v + r(s,y),\\[10pt]
v(y,0)=\varepsilon f(y),\qquad w(y,0)=\varepsilon g(y),
\end{cases}
\end{equation}
where $r$ is given by
\[
r(s,y)=\frac{b'(t)}{b(t)^2}\,w + e^{\frac N\sigma(1+\frac{\sigma}{N}-p)s}\,|v|^p.
\]

Let $\alpha(s)=\int_{\R^N} v(y,s)\,dy$. Applying integration by parts, with \eqref{system-vw}, we have
\begin{equation}\label{alpha-eq}
\begin{cases}
\displaystyle\frac{d\alpha}{ds} = \int_{\R^N} w\,dy,\\[12pt]
\displaystyle\frac{e^{-s}}{b(t)^2}\frac{d^2\alpha}{ds^2}
= \frac{e^{-s}}{b(t)^2}\frac{d\alpha}{ds} - \frac{d\alpha}{ds} + \int_{\R^N} r(s,y)\,dy .
\end{cases}
\end{equation}

Define the function $\varphi_0$ and $\psi_0$ through their Fourier transforms
\[
\widehat{\varphi}_0(\xi)=e^{-|\xi|^\sigma},\qquad
\widehat{\psi}_0(\xi)=-|\xi|^\sigma e^{-|\xi|^\sigma}.
\]
They satisfy the identities
\begin{align}
(-\Delta)^{\sigma/2}\varphi_0(y) &= \frac{y}{\sigma}\cdot\nabla_y\varphi_0(y) + \frac N\sigma\,\varphi_0(y), \label{id1}\\
\int_{\R^N}\varphi_0(y)\,dy &= 1,\qquad
\int_{\R^N}\psi_0(y)\,dy = 0. \label{id2}
\end{align}

The first identity~\eqref{id1} follows from the Fourier representation:
\[
\mathcal{F}\bigl[(-\Delta)^{\sigma/2}\varphi_0\bigr](\xi)=|\xi|^\sigma e^{-|\xi|^\sigma},
\]
while
\[
\mathcal{F}\Bigl[\frac{y}{\sigma}\cdot\nabla_y\varphi_0+\frac N\sigma\varphi_0\Bigr](\xi)
= \frac1\sigma\,i\xi\cdot i\nabla_\xi e^{-|\xi|^\sigma}
= -\frac1\sigma\,\xi\cdot\bigl(-\sigma|\xi|^{\sigma-1}\tfrac{\xi}{|\xi|}e^{-|\xi|^\sigma}\bigr)
= |\xi|^\sigma e^{-|\xi|^\sigma}.
\]
And the second identity~\eqref{id2} is obtained by evaluating the Fourier transform at $\xi=0$:
\[
\int_{\R^N}\varphi_0(y)\,dy = \widehat{\varphi}_0(0)=1,\qquad
\int_{\R^N}\psi_0(y)\,dy = \widehat{\psi}_0(0)=0.
\]

We decompose the solution $(v,w)$ as
\begin{equation}\label{decomp}
\begin{cases}
v(y,s) = \alpha(s)\,\varphi_0(y) + f(y,s),\\[4pt]
w(y,s) = \displaystyle\frac{d\alpha}{ds}(s)\,\varphi_0(y) + \alpha(s)\,\psi_0(y) + g(y,s).
\end{cases}
\end{equation}

Substituting \eqref{decomp} into \eqref{system-vw} and applying \eqref{id1}, we derive the following equations for $f$ and $g$:
\begin{equation}\label{fg-system}
\begin{cases}
f_s - \displaystyle\frac{y}{\sigma}\cdot\nabla_y f - \frac N\sigma f = g,\\[10pt]
\displaystyle\frac{e^{-s}}{b(t)^2}\Bigl(g_s - \frac{y}{\sigma}\cdot\nabla_y g - \Bigl(\frac N\sigma+1\Bigr)g\Bigr) + g = -(-\Delta_y)^{\sigma/2}f + h,
\end{cases}
\end{equation}
where the inhomogeneity $h$ is defined by:
\begin{equation}\label{h-def}
\begin{aligned}
h(y,s) &= \frac{e^{-s}}{b(t)^2}\Bigl[ -2\frac{d\alpha}{ds}\psi_0 + \alpha(s)\Bigl(\frac{y}{\sigma}\cdot\nabla_y\psi_0+\Bigl(\frac N\sigma+1\Bigr)\psi_0\Bigr)\Bigr] \\
&\quad + r(y,s) - \varphi_0(y)\int_{\R^N} r(y,s)\,dy .
\end{aligned}
\end{equation}
Integrating both sides of \eqref{h-def} with respect to $y$, and using property \eqref{id2}, together with integration by parts, we obtain
\begin{equation}\label{mean_zero_ of_h}
\int_{\R^N} h(s,y)\,dy = 0.
\end{equation}
Moreover, by construction,
\begin{equation}\label{mean_zero}
\int_{\R^N} f\,dy = \int_{\R^N} g\,dy = 0.
\end{equation}

In what follows we study the energy of the reduced system \eqref{fg-system} separately for $N=1$ and $N\ge2$.  Because the zero‑mean condition allows us to ``integrate'' the unknowns, we introduce primitives in space for $N=1$ and a Fourier‑weighted version for $N\ge2$.

\subsection{The one‑dimensional case $N=1$}

For $N=1$, we define
\[
F(y,s)=\int_{-\infty}^y f(z,s)\,dz,\qquad
G(y,s)=\int_{-\infty}^y g(z,s)\,dz,\qquad
H(y,s)=\int_{-\infty}^y h(z,s)\,dz .
\]
Then, it follows from \eqref{fg-system} that $F,G,H$ satisfy
\begin{equation}\label{FG-system}
\begin{cases}
F_s - \dfrac{y}{\sigma}F_y = G,\\[10pt]
\dfrac{e^{-s}}{b(t)^2}\bigl(G_s - \dfrac{y}{\sigma}G_y - G\bigr) + G = -(-\Delta)^{\sigma/2}F + H .
\end{cases}
\end{equation}

We now introduce three energies:
\begin{align}
E_0(s) &= \int_{\R} \Bigl( \frac12 \bigl((-\Delta)^{\sigma/4}F\bigr)^2 + \frac12\frac{e^{-s}}{b(t)^2} G^2 + \frac12 F^2 + \frac{e^{-s}}{b(t)^2} FG \Bigr)\,dy, \label{eq:E0_1d}\\
E_1(s) &= \int_{\R} \Bigl( \frac12\bigl((-\Delta)^{\sigma/4}f\bigr)^2 + \frac12\frac{e^{-s}}{b(t)^2} g^2 \Bigr) 
          + \Bigl(\frac1{2\sigma}+1\Bigr)\Bigl( \frac12 f^2 + \frac{e^{-s}}{b(t)^2} fg \Bigr)\,dy, \label{eq:E1_1d}\\
E_2(s) &= \int_{\R} y^2 \Bigl( \frac12 \bigl((-\Delta)^{\sigma/4}f\bigr)^2 + \frac12\frac{e^{-s}}{b(t)^2} g^2 + \frac12 f^2 + \frac{e^{-s}}{b(t)^2} fg \Bigr)\,dy. \label{eq:E2_1d}
\end{align}

We differentiate each energy with respect to $s$ and use the equations \eqref{FG-system} for $E_0$ or the system \eqref{fg-system} for $f,g$ for $E_1,E_2$.
All spatial integrals involving the scaling operator $y\partial_y$ are handled by integration by parts in $y$, while the terms containing the fractional Laplacian
$(-\Delta)^{\sigma/2}$ are transformed to Fourier space via Plancherel's theorem. We thus derive the following differential equations.

\begin{theorem}[Energy identities for $N=1$]\label{thm:1d_energies}
The energies defined in \eqref{eq:E0_1d}--\eqref{eq:E2_1d} satisfy
\begin{equation}
\frac{d}{ds}E_i + \frac1\sigma E_i + L_i = R_i, \label{eq:dEi_1d}
\end{equation}
for $i=0,1,2$, where the dissipative terms are
\begin{align}
L_0(s) &= \frac12\int_{\R} \bigl((-\Delta)^{\sigma/4}F\bigr)^2\,dy + \int_{\R} G^2\,dy,\\
L_1(s) &= \frac32\int_{\R} \bigl((-\Delta)^{\sigma/4}f\bigr)^2\,dy + \int_{\R} g^2\,dy 
          - \frac1\sigma\Bigl(\frac1{2\sigma}+2\Bigr)\int_{\R} f^2\,dy,\\
L_2(s) &= \int_{\R} y^2\Bigl(\frac12\bigl((-\Delta)^{\sigma/4}f\bigr)^2 + g^2\Bigr)\,dy 
          + \sigma\int_{\R} y(f+g)(-\Delta)^{\sigma/2-1}\partial_y f\,dy \nonumber\\
&\qquad - \frac{\sigma}2\Bigl(\frac{\sigma}2-1\Bigr)\int_{\R} (f+g)(-\Delta)^{\sigma/2-1}f\,dy,
\end{align}
and the remainders are given explicitly by
\begin{align}
R_0(s) &= \frac{3}{2}\frac{e^{-s}}{b(t)^2}\int_{\R} G^2\,dy 
        - \frac{2}{b(t)^2}\frac{db}{dt}\int_{\R} FG\,dy 
        - \frac{1}{b(t)^2}\frac{db}{dt}\int_{\R} G^2\,dy 
        + \int_{\R} (F+G)H\,dy,\\[4pt]
R_1(s) &= \Bigl(\frac52+\frac{2}{\sigma}\Bigr)\frac{e^{-s}}{b(t)^2}\int_{\R} g^2\,dy 
        + \frac{2}{\sigma}\Bigl(\frac1{2\sigma}+1\Bigr)\frac{e^{-s}}{b(t)^2}\int_{\R} fg\,dy \nonumber\\
&\quad - \frac{1}{b(t)^2}\frac{db}{dt}\int_{\R} \Bigl(4\Bigl(\frac1{2\sigma}+1\Bigr)fg + g^2\Bigr)\,dy 
        + \int_{\R} \Bigl(2\Bigl(\frac1{2\sigma}+1\Bigr)f + g\Bigr)h\,dy,\\[4pt]
R_2(s) &= \frac{3}{2}\frac{e^{-s}}{b(t)^2}\int_{\R} y^2 g^2\,dy 
        - \frac{1}{b(t)^2}\frac{db}{dt}\int_{\R} y^2(2f+g)g\,dy 
        + \int_{\R} y^2(f+g)h\,dy.
\end{align}
\end{theorem}

\subsection{The higher‑dimensional case $N\ge2$}

Choose a small parameter $\delta\in(0,\sigma)$ and set
\[
\widehat{F}(s,\xi)=|\xi|^{-\frac N2-\delta}\,\widehat{f}(s,\xi),\qquad
\widehat{G}(s,\xi)=|\xi|^{-\frac N2-\delta}\,\widehat{g}(s,\xi),\qquad
\widehat{H}(s,\xi)=|\xi|^{-\frac N2-\delta}\,\widehat{h}(s,\xi).
\]
Then the system for $(\widehat F,\widehat G)$ becomes
\begin{equation}\label{eq:FG_ND}
\begin{cases}
\widehat F_s + \displaystyle\frac{\xi}{\sigma}\cdot\nabla_\xi\widehat F
+ \frac1\sigma\Bigl(\frac N2+\delta\Bigr)\widehat F = \widehat G,\\[12pt]
\displaystyle\frac{e^{-s}}{b(t)^2}\Bigl(\widehat G_s + \frac{\xi}{\sigma}\cdot\nabla_\xi\widehat G
+ \frac1\sigma\Bigl(\frac N2+\delta-\sigma\Bigr)\widehat G\Bigr) + \widehat G
= -|\xi|^\sigma\widehat F + \widehat H .
\end{cases}
\end{equation}

We now introduce three energies:
\begin{align}
E_0(s) &= \operatorname{Re}\int_{\R^N} \Bigl( \frac12\bigl(|\xi|^\sigma|\widehat{F}|^2 + \frac{e^{-s}}{b^2}|\widehat{G}|^2\bigr) 
                  + \frac12|\widehat{F}|^2 + \frac{e^{-s}}{b^2}\widehat{F}\overline{\widehat{G}} \Bigr)\,d\xi, \label{eq:calE0}\\
E_1(s) &= \int_{\R^N} \Bigl( \frac12\bigl(|\xi|^\sigma|\widehat{F}|^2 + \frac{e^{-s}}{b^2}|\widehat{G}|^2\bigr) \nonumber\\
&\qquad\qquad + \Bigl(\frac N{2\sigma}+1\Bigr)\Bigl(\frac12|\widehat{F}|^2 + \frac{e^{-s}}{b^2}\operatorname{Re}(\widehat{F}\overline{\widehat{G}})\Bigr) \Bigr)\,d\xi, \label{eq:calE1}\\
E_2(s) &= \int_{\R^N} |y|^{2m} \Bigl[ \frac12 \Bigl( ((-\Delta)^{\sigma/4}f)^2 + \frac{e^{-s}}{b(t)^2} g^2 \Bigr) 
          + \frac12 f^2 + \frac{e^{-s}}{b(t)^2} fg \Bigr]\,dy. \label{eq:E2_ND}
\end{align}

Similarly to the one‑dimensional case, we differentiate $E_0,E_1$ in Fourier space and $E_2$ in physical space, using respectively the system \eqref{eq:FG_ND} and the equations for $f,g$. All terms involving the scaling operator $\xi\cdot\nabla_\xi$ or $y\cdot\nabla_y$ are handled by integration by parts, while the fractional Laplacian is treated through the multiplier $|\xi|^\sigma$. We thus obtain the following theorem.

\begin{theorem}[Energy identities for $N\ge2$]\label{thm:ND_energies}
The energies defined in \eqref{eq:calE0}--\eqref{eq:E2_ND} satisfy
\begin{align}
\frac{d}{ds}E_0 + \frac{2\delta}{\sigma}E_0 + L_0 &= R_0, \label{eq:dE0_ND}\\
\frac{d}{ds}E_1 + \frac{2\delta}{\sigma}E_1 + L_1 &= R_1, \label{eq:dE1_ND}\\
\frac{d}{ds}E_2 + \frac{2m-N-\kappa}{\sigma}E_2 + L_2 &= R_2, \label{eq:dE2_ND}
\end{align}
where $\kappa>0$ is a free parameter,
\begin{align}
L_0(s) &= \frac12\int_{\R^N} |\xi|^\sigma|\widehat{F}|^2\,d\xi + \int_{\R^N} |\widehat{G}|^2\,d\xi,\\
L_1(s) &= \frac12\Bigl(1-\frac{2\delta}{\sigma}\Bigr)\int_{\R^N} |\xi|^\sigma|\widehat{F}|^2\,d\xi 
                  + \int_{\R^N} |\widehat{G}|^2\,d\xi \nonumber\\
&\quad - \frac12\Bigl(\frac N\sigma+\frac{2\delta}{\sigma}\Bigr)\Bigl(\frac N{2\sigma}+1\Bigr)\int_{\R^N} |\widehat{F}|^2\,d\xi,\\
L_2(s) &= \frac{\kappa}{2\sigma}\int_{\R^N}|y|^{2m}f^2\,dy
        + \frac{\sigma+\kappa}{2\sigma}\int_{\R^N}|y|^{2m}\bigl((-\Delta)^{\sigma/4}f\bigr)^2\,dy
        - \int_{\R^N}|y|^{2m}\bigl((-\Delta)^{\sigma/4}f\bigr)^2\,dy \nonumber\\
&\quad + \int_{\R^N}|y|^{2m}(-\Delta)^{\sigma/2}f\cdot f\,dy
        + \int_{\R^N}|y|^{2m}g^2\,dy
        + \int_{\R^N}|y|^{2m}g\,(-\Delta)^{\sigma/2}f\,dy\\
&\quad - \int_{\R^N}|y|^{2m}(-\Delta)^{\sigma/4}f\,(-\Delta)^{\sigma/4}g\,dy,\nonumber
\end{align}
and the remainders are given explicitly by
\begin{align}
R_0(s) &= \frac{3}{2}\frac{e^{-s}}{b(t)^2}\int_{\R^N} |\widehat{G}|^2\,d\xi 
        - \frac{2}{b(t)^2}\frac{db}{dt}\operatorname{Re}\int_{\R^N} \widehat{F}\overline{\widehat{G}}\,d\xi 
        - \frac{1}{b(t)^2}\frac{db}{dt}\int_{\R^N} |\widehat{G}|^2\,d\xi \nonumber\\
&\quad + \operatorname{Re}\int_{\R^N} (\widehat{F}+\widehat{G})\overline{\widehat{H}}\,d\xi,\\[4pt]
R_1(s) &= \Bigl(\frac52+\frac{2}{\sigma}\Bigr)\frac{e^{-s}}{b(t)^2}\int_{\R^N} |\widehat{G}|^2\,d\xi 
        + \frac{2}{\sigma}\Bigl(\frac N{2\sigma}+1\Bigr)\frac{e^{-s}}{b(t)^2}\operatorname{Re}\int_{\R^N} \widehat{F}\overline{\widehat{G}}\,d\xi \nonumber\\
&\quad - \frac{1}{b(t)^2}\frac{db}{dt}\int_{\R^N} \Bigl(4\Bigl(\frac N{2\sigma}+1\Bigr)\operatorname{Re}(\widehat{F}\overline{\widehat{G}}) + |\widehat{G}|^2\Bigr)\,d\xi\\
&\quad + \operatorname{Re}\int_{\R^N} \Bigl(2\Bigl(\frac N{2\sigma}+1\Bigr)\widehat{F} + \widehat{G}\Bigr)\overline{\widehat{H}}\,d\xi, \nonumber\\[4pt]
R_2(s) &= \frac{3}{2}\frac{e^{-s}}{b(t)^2}\int_{\R^N} |y|^{2m} g^2\,dy 
        - \frac{1}{b(t)^2}\frac{db}{dt}\int_{\R^N} |y|^{2m}(2f+g)g\,dy 
        + \int_{\R^N} |y|^{2m}(f+g)h\,dy.
\end{align}
\end{theorem}

Moreover, by using Lemma 3.2 of \cite{Wakasugi2017}, the following equivalents are valid for $s \ge s_1$ with sufficiently large $s_1 > 0$. When $N=1$,
\begin{align*}
    E_0(s) &\sim \int_{\mathbb{R}} \left( \frac{1}{2} ((-\Delta)^{\frac{\sigma}{4}} F)^2 + \frac{e^{-s}}{b(t)^2} G^2 \right) + \frac{1}{2} F^2 \, dy, \\[1em]
    E_1(s) &\sim \int_{\mathbb{R}} \left[ \frac{1}{2} \left( ((-\Delta)^{\frac{\sigma}{4}} f)^2 + \frac{e^{-s}}{b(t)^2} g^2 \right) + \left( \frac{1}{2\sigma} + 1 \right)  f^2 \right] dy, \\[1em]
    E_2(s) &\sim \int_{\mathbb{R}} y^2 \left[ \left( \frac{1}{2} ((-\Delta)^{\frac{\sigma}{4}} f)^2 + \frac{e^{-s}}{b(t)^2} g^2 \right) + \frac{1}{2} f^2 \right] dy.
\end{align*}
When $N\ge 2$,
\begin{align*}
E_0(s) &\sim \int_{\mathbb{R}^N} \left( \frac{1}{2} \left( |\xi|^\sigma |\widehat{F}|^2 + \frac{e^{-s}}{b(t)^2} |\widehat{G}|^2 \right) + \frac{1}{2} |\widehat{F}|^2\right) d\xi, \\[1em]
E_1(s) &\sim \int_{\mathbb{R}^N}  \frac{1}{2}\left( \left( (-\Delta)^{\frac{\sigma}{4}} f \right)^2 + \frac{e^{-s}}{b(t)^2} g^2 \right) + \left( \frac{N}{4\sigma} + 1 \right)  f^2 \, dy, \\[1em]
E_2(s) &\sim \int_{\mathbb{R}^N} |y|^{2m} \left[ \frac{1}{2} \left( ((-\Delta)^{\frac{\sigma}{4}} f)^2 + \frac{e^{-s}}{b(t)^2} g^2 \right) + \frac{1}{2} f^2 \right] dy.
\end{align*}

\subsection{A prior estimate}

For a parameter $\lambda>0$, define
\[
E_3(s) = \frac12\frac{e^{-s}}{b(t)^2}\Bigl(\frac{d\alpha}{ds}\Bigr)^2 + e^{-2\lambda s}\alpha(s)^2 .
\]
A direct differentiation yields
\begin{equation}\label{E3-derivative}
\frac{d}{ds}E_3(s) + 2\lambda E_3(s) + \Bigl(\frac{d\alpha}{ds}\Bigr)^2 = R_3(s),
\end{equation}
with
\[
R_3(s) = \frac12(2\lambda+1)\frac{e^{-s}}{b(t)^2}\Bigl(\frac{d\alpha}{ds}\Bigr)^2
- \frac1{b(t)^2}\frac{db}{dt}\Bigl(\frac{d\alpha}{ds}\Bigr)^2
+ \frac{d\alpha}{ds}\int_{\R^N} r\,dy
+ e^{-2\lambda s}\alpha\frac{d\alpha}{ds}.
\]

Then, choose positive constants
$C_0\gg C_1\gg 1$ and set
\begin{equation}\label{E4}
E_4(s) = C_0E_0(s) + C_1E_1(s) + E_2(s) + E_3(s).
\end{equation}
From the differential identities in Theorem \ref{thm:1d_energies} and Theorem \ref{thm:ND_energies}, we obtain an equation of the form
\begin{equation}\label{eq:E4}
\frac{d}{ds}E_4(s) + 2\lambda E_4(s) + L_4(s) = R_4(s),
\end{equation}
where 
\[
\begin{aligned}
L_4(s) &= (\frac{1}{\sigma}-2\lambda)(C_0E_0(s) + C_1E_1(s) + E_2(s) + E_3(s))\\
&\quad + C_0L_0(s) + C_1L_1(s) + L_2(s) + \left(\frac{d\alpha}{ds}\right)^2
\end{aligned}
\]
for $N=1$, and
\[
\begin{aligned}
L_4(s)&=C_0\Bigl(\frac{2\delta}{\sigma}-2\lambda\Bigr)E_0(s)
+ C_1\Bigl(\frac{2\delta}{\sigma}-2\lambda\Bigr)E_1(s)
+ \Bigl(\frac1\sigma(2m-N-\kappa)-2\lambda\Bigr)E_2(s)\\
&\quad+ C_0L_0(s)+C_1L_1(s)+L_2(s)+\Bigl(\frac{d\alpha}{ds}\Bigr)^2
\end{aligned}
\]
for $N\ge 2$, and
\[
R_4 = C_0R_0(s) + C_1R_1(s) + R_2(s) + R_3(s).
\]
By choosing $\lambda$ sufficiently small and the parameters $\delta,\kappa,m$ appropriately, one can make $L_4$ strictly positive and absorb the remainder $R_4$ into $L_4$ and $E_4$, provided that the damping coefficient $b(t)$ enjoys suitable asymptotic properties.

Noticing the asymptotic equivalence of $E_0, E_1$, and $E_2$, it suffices to estimate the non-positive components within $L_0, L_1$, and $L_2$. 

First, we consider the case \(N=1\). To estimate the non-positive term \(\int_{\mathbb{R}} f^2\,dy\), we introduce the following lemma.

\begin{lemma}
\label{lem:low_freq}
For the one-dimensional function $f$, the following estimate holds:
\begin{equation}
\int_{\mathbb{R}} f^2 dy \le \eta \int_{\mathbb{R}} |\xi|^\sigma |\widehat{f}|^2 d\xi + \eta^{1 - \frac{2}{\sigma}} \int_{\mathbb{R}} |(-\Delta)^{\frac{\sigma}{4}}F|^2 dy
\end{equation}
\end{lemma}

\begin{proof}
By Plancherel's theorem and a frequency splitting, we have
\begin{equation}
\int_{\mathbb{R}} f^2 dy = \int_{\mathbb{R}} |\widehat{f}|^2 d\xi = \int_{|\xi| \ge \eta^{-\frac{1}{\sigma}}} |\widehat{f}|^2 d\xi + \int_{|\xi| < \eta^{-\frac{1}{\sigma}}} |\widehat{f}|^2 d\xi
\end{equation}
For the high-frequency, we directly extract the differential weight:
\begin{equation}
\int_{|\xi| \ge \eta^{-\frac{1}{\sigma}}} |\widehat{f}|^2 d\xi \le \left(\eta^{-\frac{1}{\sigma}}\right)^{-\sigma} \int_{|\xi| \ge \eta^{-\frac{1}{\sigma}}} |\xi|^\sigma |\widehat{f}|^2 d\xi \le \eta \int_{\mathbb{R}} |\xi|^\sigma |\widehat{f}|^2 d\xi
\end{equation}
For the low-frequency regime ($|\xi| < \eta^{-\frac{1}{\sigma}}$), applying the differential relation $i\xi \widehat{F}(\xi) = \widehat{f}(\xi)$, we deduce
\begin{equation}
\int_{|\xi| < \eta^{-\frac{1}{\sigma}}} |\widehat{f}|^2 d\xi \le \left(\eta^{-\frac{1}{\sigma}}\right)^{2-\sigma} \int_{|\xi| < \eta^{-\frac{1}{\sigma}}} |\xi|^\sigma |\widehat{F}|^2 d\xi = \eta^{1-\frac{2}{\sigma}} \int_{\mathbb{R}} |(-\Delta)^{\frac{\sigma}{4}}F|^2 dy
\end{equation}
Combining the estimates from both the high-frequency and low-frequency regimes directly yields the desired inequality, which completes the proof.
\end{proof}

With Lemma~\ref{lem:low_freq} in hand, we now turn to the non-positive terms appearing in the \(L_2\). Their control is provided by the next lemma.
\begin{lemma}
\label{lem:nonlocal_cross}
By splitting the frequency domain at \(|\xi|=\eta\), the non-positive terms involving \(g\) and \(f\) satisfy the following bounds:
\begin{align}
\left| \int_{\mathbb{R}} y g (-\Delta)^{\frac{\sigma}{2}-1} \partial_y f dy \right| &\le \frac{1}{4\sigma}\eta^\sigma \|yf\|_{L^2}^2 + \frac{1}{2\sigma}\eta^\sigma \|yg\|_{L^2}^2 + \frac{1}{2}\sigma \eta^{-2} \left\|(-\Delta)^{\frac{\sigma}{2}}f\right\|_{L^2}^2 \\
\left| \int_{\mathbb{R}} g (-\Delta)^{\frac{\sigma-1}{2}} f dy \right| &\le \frac{1}{2\sigma}\eta^\sigma \|yg\|_{L^2}^2 + \frac{1}{2\sigma}\eta^\sigma \|yf\|_{L^2}^2 + \frac{1}{2}\eta^{\frac{\sigma}{2}-2} \|g\|_{L^2}^2 + \frac{1}{2}\eta^{\frac{\sigma}{2}-2} \left\|(-\Delta)^{\frac{\sigma}{4}}f\right\|_{L^2}^2
\end{align}
\end{lemma}

\begin{proof}
\textbf{Part I. Estimate of $\int_{\mathbb{R}} y g (-\Delta)^{\frac{\sigma}{2}-1} \partial_y f \, dy$}

Applying the Cauchy--Schwarz inequality yields
\begin{equation}
\int_{\mathbb{R}} y g (-\Delta)^{\frac{\sigma}{2}-1} \partial_y f \, dy \le \|y g\|_{L^2} \left\| (-\Delta)^{\frac{\sigma}{2}-1} \partial_y f \right\|_{L^2}
\end{equation}

By Plancherel's theorem, the second norm can be rewritten and split in frequency space as
\begin{equation}
\left\| (-\Delta)^{\frac{\sigma}{2}-1} \partial_y f \right\|_{L^2}^2 = \int_{\mathbb{R}} |\xi|^{2\sigma-2} |\hat{f}|^2 \, d\xi = \int_{|\xi| \le \eta} |\xi|^{2\sigma-2} |\hat{f}|^2 \, d\xi + \int_{|\xi| > \eta} |\xi|^{2\sigma-2} |\hat{f}|^2 \, d\xi
\end{equation}

For \(|\xi|>\eta\), we have
\begin{equation}
\begin{aligned}
\int_{|\xi| > \eta} |\xi|^{2\sigma-2} |\hat{f}|^2 \, d\xi &< \eta^{\sigma-2} \int_{|\xi| > \eta} |\xi|^{\sigma} |\hat{f}|^2 \, d\xi \nonumber \\
&< \eta^{\sigma-2} \left\| (-\Delta)^{\frac{\sigma}{4}} f \right\|_{L^2}^2.
\end{aligned}
\end{equation}

For \(|\xi|\le\eta\), we use the zero-mean condition \(\hat{f}(0)=\int_{\mathbb{R}} f(x)\,dx=0\) to write 
\[
\hat{f}(\xi)=\int_0^\xi \hat{f}'(z)\,dz,
\]
which implies
\begin{equation}
|\hat{f}(\xi)| \le |\xi|^{1/2}\|\hat{f}'\|_{L^2} = |\xi|^{1/2}\|yf\|_{L^2}.
\end{equation}
Thus,
\begin{equation}\label{lowf}
\begin{aligned}
\int_{|\xi|\le\eta} |\xi|^{2\sigma-2}|\hat{f}|^2\,d\xi
&\le \int_{|\xi|\le\eta} |\xi|^{2\sigma-1}\|yf\|_{L^2}^2\,d\xi \\
&= \frac{1}{2\sigma}\eta^{2\sigma}\|yf\|_{L^2}^2 .
\end{aligned}
\end{equation}
Combining these estimates, we obtain

\begin{equation}
\left\| (-\Delta)^{\frac{\sigma}{2}-1} \partial_y f \right\|_{L^2}^2 \le \frac{1}{2\sigma} \eta^{2\sigma} \|y f\|_{L^2}^2 + \eta^{\sigma-2} \left\| (-\Delta)^{\frac{\sigma}{4}} f \right\|_{L^2}^2
\end{equation}
i.e.,
\begin{equation}
\left\| (-\Delta)^{\frac{\sigma}{2}-1} \partial_y f \right\|_{L^2} \le \frac{1}{2\sqrt{\sigma}} \eta^\sigma \|y f\|_{L^2} + \frac{\sqrt{2}}{2}\eta^{\frac{\sigma}{2}-1} \left\| (-\Delta)^{\frac{\sigma}{2}} f \right\|_{L^2}.
\end{equation}

Consequently, by Young's inequality, we obtain
\begin{equation}
\begin{aligned}
&\int_{\mathbb{R}} y g (-\Delta)^{\frac{\sigma}{2}-1} \partial_y f \, dy\\
&\le \frac{1}{2\sqrt{\sigma}} \eta^\sigma \|y f\|_{L^2} \|y g\|_{L^2} + \frac{\sqrt{2}}{2}\eta^{\frac{\sigma}{2}-1} \left\| (-\Delta)^{\frac{\sigma}{2}} f \right\|_{L^2} \|y g\|_{L^2} \\
&\le \frac{1}{4\sigma} \eta^\sigma \|y f\|_{L^2}^2 + \frac{1}{4\sigma} \eta^\sigma \|y g\|_{L^2}^2 + \frac{1}{4\sigma} \eta^\sigma \|y g\|_{L^2}^2 + \frac{1}{2} \sigma \eta^{-2} \left\| (-\Delta)^{\frac{\sigma}{2}} f \right\|_{L^2}^2 \\
&= \frac{1}{4\sigma} \eta^\sigma \|y f\|_{L^2}^2 + \frac{1}{2\sigma} \eta^\sigma \|y g\|_{L^2}^2 + \frac{1}{2} \sigma \eta^{-2} \left\| (-\Delta)^{\frac{\sigma}{2}} f \right\|_{L^2}^2
\end{aligned}
\end{equation}

\textbf{Part II. Estimate of $\int_{\mathbb{R}} g (-\Delta)^{\frac{\sigma}{2}-1} f \, dy$}

Using Plancherel's theorem,
\begin{align}
\int_{\mathbb{R}} g (-\Delta)^{\frac{\sigma}{2}-1} f \, dy &= \int_{\mathbb{R}} \bar{\hat{g}} \, |\xi|^{\sigma-2} \, \hat{f} \, d\xi \nonumber \\
&\le \int_{\mathbb{R}} |\xi|^{\sigma-2} |\hat{g}| |\hat{f}| \, d\xi \\
&= \int_{|\xi| \le \eta} |\xi|^{\sigma-2} |\hat{g}| |\hat{f}| \, d\xi + \int_{|\xi| > \eta} |\xi|^{\sigma-2} |\hat{g}| |\hat{f}| \, d\xi \nonumber 
\end{align}

For \(|\xi|\le\eta\), the Cauchy--Schwarz inequality together with an estimate similar to \eqref{lowf} gives 
\begin{equation}
\begin{aligned}
\int_{|\xi| \le \eta} |\xi|^{\sigma-2} |\hat{g}| |\hat{f}| \, d\xi &\le \left( \int_{|\xi| \le \eta} |\xi|^{\sigma-2} |\hat{g}|^2 \, d\xi \right)^{1/2} \left( \int_{|\xi| \le \eta} |\xi|^{\sigma-2} |\hat{f}|^2 \, d\xi \right)^{1/2} \\
&\le \frac{1}{2\sigma} \eta^\sigma \|y g\|_{L^2} \|y f\|_{L^2}
\end{aligned}
\end{equation}

For \(|\xi|>\eta\), shifting the frequency weights yields
\begin{equation}
\begin{aligned}
\int_{|\xi| > \eta} |\xi|^{\sigma-2} |\hat{g}| |\hat{f}| \, d\xi &< \eta^{\frac{\sigma}{2}-2} \int_{|\xi| > \eta} |\xi|^{\frac{\sigma}{2}} |\hat{g}| |\hat{f}| \, d\xi \\
&< \eta^{\frac{\sigma}{2}-2} \|g\|_{L^2} \left\| (-\Delta)^{\frac{\sigma}{4}} f \right\|_{L^2}
\end{aligned}
\end{equation}

Thus, putting the two regimes together, we have
\begin{align}
\int_{\mathbb{R}} g (-\Delta)^{\frac{\sigma}{2}-1} f \, dy &< \frac{1}{2\sigma} \eta^\sigma \|y g\|_{L^2} \|y f\|_{L^2} + \eta^{\frac{\sigma}{2}-2} \|g\|_{L^2} \left\| (-\Delta)^{\frac{\sigma}{4}} f \right\|_{L^2} \nonumber \\
&\le \frac{1}{2\sigma} \eta^\sigma \|y g\|_{L^2}^2 + \frac{1}{2\sigma} \eta^\sigma \|y f\|_{L^2}^2 \\
&\quad + \frac{1}{2} \eta^{\frac{\sigma}{2}-2} \|g\|_{L^2}^2 + \frac{1}{2} \eta^{\frac{\sigma}{2}-2} \left\| (-\Delta)^{\frac{\sigma}{4}} f \right\|_{L^2}^2 \nonumber 
\end{align}
\end{proof}

Next, we consider the higher-dimension case \(N \ge 2\). The same non-positive term \(\int_{\mathbb{R}} f^2\,dy\) could be estimated as follows.

\begin{lemma}
\label{lem:low_freq_N}
For the function $f$ in higher-dimension case \(N \ge 2\), the following estimate holds:
\begin{equation}
\int_{\mathbb{R}} f^2 dy \le \eta \int_{\mathbb{R}} |\xi|^\sigma |\widehat{f}|^2 d\xi + \eta^{1 - \frac{2}{\sigma}} \int_{\mathbb{R}} |(-\Delta)^{\frac{\sigma}{4}}F|^2 dy
\end{equation}
\end{lemma}

\begin{proof}
Proceeding as in the proof of Lemma~\ref{lem:low_freq}, we deduce
\[
\begin{aligned}
\int_{\mathbb{R}} f^2 dy &= \int_{\mathbb{R}} |\widehat{f}|^2 d\xi\\
& = \int_{|\xi| \ge \eta^{-\frac{1}{\sigma}}} |\widehat{f}|^2 d\xi + \int_{|\xi| < \eta^{-\frac{1}{\sigma}}} |\widehat{f}|^2 d\xi\\
& \le \eta \int_{\mathbb{R}} |\xi|^\sigma |\widehat{f}|^2 d\xi + \eta^{\frac{\sigma-N-2\delta}{\sigma}} \int_{|\xi| < \eta^{-\frac{1}{\sigma}}} |\xi|^{\sigma-N-2\delta} |\widehat{f}|^2 d\xi\\
& \le \eta \int_{\mathbb{R}} |\xi|^\sigma |\widehat{f}|^2 d\xi + \eta^{\frac{\sigma-N-2\delta}{\sigma}} \int_{|\xi| < \eta^{-\frac{1}{\sigma}}} |\xi|^{\sigma} |\widehat{F}|^2 d\xi,
\end{aligned}
\]
since \(\sigma-N-2\delta<0\) for \(N\ge 2\).
\end{proof}

As in the one-dimensional case, we have to estimate the following two differences:
\[
\int_{\R^N}|y|^{2m}(-\Delta)^{\sigma/2}f\cdot f\,dy - \int_{\R^N}|y|^{2m}\bigl((-\Delta)^{\sigma/4}f\bigr)^2\,dy
\]
and
\[
\int_{\R^N}|y|^{2m}g\,(-\Delta)^{\sigma/2}f\,dy - \int_{\R^N}|y|^{2m}(-\Delta)^{\sigma/4}f\,(-\Delta)^{\sigma/4}g\,dy.
\]

Applying Lemma~\ref{thm:weighted_commutator} with $u(x)=(-\Delta)^{\sigma/4}f$ and $U(x)=f$, we get
\begin{equation}
\begin{aligned}
&\left|\int_{\R^N}|y|^{2m}(-\Delta)^{\sigma/2}f\cdot f\,dy - \int_{\R^N}|y|^{2m}\bigl((-\Delta)^{\sigma/4}f\bigr)^2\,dy\right|\\
& = \left|\int_{\R^N}\left([(-\Delta)^{\sigma/4},|y|^{2m}](-\Delta)^{\sigma/4}f\right) f\,dy \right|\\
&\le
C \norm{\la x\ra^{m-\varepsilon}(-\Delta)^{\sigma/4}f}_{L^2} \norm{\la x\ra^m f}_{L^2}
+ C \norm{\la x\ra^m (-\Delta)^{\sigma/4}f}_{L^2} \norm{\la x\ra^{m-\varepsilon}f}_{L^2}\\
& \le \eta \norm{\la x\ra^{m}(-\Delta)^{\sigma/4}f}_{L^2}^2+C_{\eta}\norm{(-\Delta)^{\sigma/4}f}_{L^2}^2+C \norm{\la x\ra^m f}_{L^2}^2.
\end{aligned}
\end{equation}

For the second term, we have the similar estimate
\begin{equation}
\begin{aligned}
&\left|\int_{\R^N}|y|^{2m}g\,(-\Delta)^{\sigma/2}f\,dy - \int_{\R^N}|y|^{2m}(-\Delta)^{\sigma/4}f\,(-\Delta)^{\sigma/4}g\,dy\right|\\
& \le \eta \norm{\la x\ra^{m}(-\Delta)^{\sigma/4}f}_{L^2}^2+C_{\eta}\norm{(-\Delta)^{\sigma/4}f}_{L^2}^2+C \norm{\la x\ra^m g}_{L^2}^2.
\end{aligned}
\end{equation}

Taking $\eta$ sufficiently small and $C_0, C_1$ sufficiently large with $C_1\ll C_0$, we could get the following estimate:
\begin{equation}
L_4(s) \ge C \left( \|f(s)\|_{H^{\sigma/2,m}}^2 + \|g(s)\|_{H^{0,m}}^2 + \left(\frac{d\alpha}{ds}(s)\right)^2 \right)
\end{equation}
for all $s > s_1$.

\subsection{Total energy}
We introduce the total energy $E_5(s)$ as:
\begin{equation}\label{eq:E5}
E_5(s) = E_4(s) + \frac{1}{2} \alpha(s)^2 + \frac{e^{-s}}{b(t)^2} \alpha(s) \frac{d\alpha}{ds}.
\end{equation}

By similar argument in \cite{Wakasugi2017}, we have
\begin{lemma}[Energy Equivalence of $E_5(s)$]\label{lem:E5_equiv}
There exists $s_2 \ge s_1$ such that for all $s \ge s_2$, the modified energy $E_5(s)$ satisfies:
\begin{equation}\label{eq:E5_equiv}
E_5(s) \sim \|f(s)\|_{H^{\sigma/2,m}}^2 + \frac{e^{-s}}{b(t)^2} \|g(s)\|_{H^{0,m}}^2 + \alpha(s)^2 + \frac{e^{-s}}{b(t)^2} \frac{d\alpha}{ds}(s)^2.
\end{equation}
\end{lemma}

Taking the derivative of $E_5(s)$ and combining it with the energy identity \eqref{eq:E4} yields
\begin{equation}\label{eq:E5_diff}
\frac{d}{ds} E_5(s) + 2\lambda E_4(s) + L_4(s) = R_5(s)
\end{equation}
where $R_5(s) = R_4(s) + \tilde{R}_5(s)$ and 
$$
\tilde{R}_5(s) = \frac{e^{-s}}{b(t)^2}\left(\frac{d\alpha}{ds}\right)^2 - \frac{2}{b(t)^2}\frac{db}{dt}\alpha(s)\frac{d\alpha}{ds} + \alpha(s)\int_{\mathbb{R}^N} r(s,y) \, dy.
$$ 

We now proceed to bound $R_5(s)$. The estimates are largely parallel to those in \cite{Wakasugi2017}, the main difference being the occurrence of the fractional Laplacian in the energies and the estimate of the nonlinear term in $H^{\sigma/2,m}$. 

First, by the Gagliardo--Nirenberg inequality in weighted fractional Sobolev spaces and Lemma 2.1 in \cite{LI}, the nonlinear term could be estimated as follows:
\begin{equation}\label{eq:nonlinear_est}
\begin{aligned}
\left\| e^{\frac{N}{\sigma}\left(1+\frac{\sigma}{N}-p\right)s} |v|^p \right\|_{H^{0,m}}^2 
&= \int_{\mathbb{R}^N} \langle y \rangle^{2m} e^{2\frac{N}{\sigma}\left(1+\frac{\sigma}{N}-p\right)s} |v|^{2p} \, dy  \\
&\le 2 e^{2\frac{N}{\sigma}\left(1+\frac{\sigma}{N}-p\right)s} \int_{\mathbb{R}^N} [\langle y \rangle^{\frac mp} |v|]^{2p} \, dy \\
&\le C e^{2\frac{N}{\sigma}\left(1+\frac{\sigma}{N}-p\right)s} \|\langle y \rangle^{\frac mp} v\|_{H^{\frac \sigma 2,0}}^{2p\theta} \|\langle y \rangle^{\frac mp} v\|_{L^2}^{2p(1-\theta)} \\
& \le C e^{2\frac{N}{\sigma}\left(1+\frac{\sigma}{N}-p\right)s} \|v\|_{H^{\frac \sigma 2,m/p}}^{2p}\\
& \le C e^{2\frac{N}{\sigma}\left(1+\frac{\sigma}{N}-p\right)s} \|v\|_{H^{\frac \sigma 2,m}}^{2p},
\end{aligned}
\end{equation}
with $\theta = \frac n\sigma\left(1-\frac 1p\right)$ and $p\leq \frac{n}{n-\sigma}$.

Using the decomposition $w = \frac{d\alpha}{ds}\varphi_0 + \alpha\psi_0 + g$ together with Lemma 3.2 in \cite{Wakasugi2017}, the coefficient terms are controlled by
\begin{equation}\label{eq:coeff_est}
\Bigl\| \frac{1}{b(t)^2} \frac{db}{dt} w \Bigr\|_{H^{0,m}} \le C e^{-\lambda_0 s} \bigl( \|g\|_{H^{0,m}} + |\alpha| + \left|\frac{d\alpha}{ds}\right| \bigr),
\end{equation}
where $\lambda_0 = \frac{1-\beta}{1+\beta}$ for $\beta\in(-1,1)$ and $\lambda_0$ can be taken arbitrarily large when $\beta=-1$.

From \eqref{eq:nonlinear_est} and \eqref{eq:coeff_est}, we obtain the bound for the inhomogeneity $r$:
\begin{equation}\label{eq:r_est}
\|r(s)\|_{H^{0,m}} \le C e^{-\lambda_0 s} \bigl( \|g\|_{H^{0,m}} + |\alpha| + |\frac{d\alpha}{ds}| \bigr) + C e^{-\lambda_1 s} \bigl( \|f\|_{H^{\sigma/2,m}} + |\alpha| \bigr)^p,
\end{equation}
where $\lambda_1 = \frac{N}{\sigma}\left( p - 1 - \frac{\sigma}{N} \right)$.

By the definition of the function $h$ and the zero-mean property $\int_{\R^N} h dy = 0$, the $H^{0,m}$-norm of $h$ satisfies the same type of estimate as $r$:
\begin{equation}\label{eq:h_est}
\|h\|_{H^{0,m}} \le C e^{-\lambda_0 s} \bigl( \|g\|_{H^{0,m}} + |\alpha| + |\alpha'| \bigr) + C e^{-\lambda_1 s} \bigl( \|f\|_{H^{\sigma/2,m}} + |\alpha| \bigr)^p .
\end{equation}
It follows from Lemma 3.9 and Lemma 3.11 of \cite{Wakasugi2017} that
\begin{equation}
\begin{aligned}
\|H(s)\|_{L^2}^2 &\le C \|h(s)\|_{H^{0,m}}^2 \\
&\le C e^{-2\lambda_0 s} \left( \|g(s)\|_{H^{0,m}}^2 + \alpha(s)^2 + \left(\frac{d\alpha}{ds}\right)^2 \right) + C e^{2\frac{N}{\sigma}\left(1+\frac{\sigma}{N}-p\right)s} \left( \|f(s)\|_{H^{\frac \sigma 2,m}} + \alpha(s) \right)^{2p}
\end{aligned}
\end{equation}

Following similar arguments as in \cite[Section~4]{Wakasugi2017}, the terms in $R_5$  can be estimated using Young's inequality together with \eqref{eq:h_est} and the equivalence of the energies. Thus, we get

\begin{equation}\label{eq:R5_estimate}
|R_5(s)| \le \tilde{\eta} L_4(s) + C(\tilde{\eta}) e^{-\lambda_0 s} E_5(s) + C(\tilde{\eta}) e^{-2\lambda_1 s} E_5(s)^p + C(\tilde{\eta}) e^{-\lambda_1 s} E_5(s)^{(p+1)/2},
\end{equation}
for any small $\tilde{\eta} > 0$.

Combining the above estimates and choosing $\tilde{\eta}$ sufficiently small to absorb the $L_4$-terms, we obtain the key differential inequality for $E_5$:
\begin{equation}\label{eq:closed_ineq}
\frac{d}{ds} E_5(s) \le C e^{-\lambda_0 s} E_5(s) + C e^{-2\lambda_1 s} E_5(s)^p + C e^{-\lambda_1 s} E_5(s)^{(p+1)/2},
\end{equation}
valid for all $s \ge s_0$ with $s_0$ large enough.

\section{Global existence in the supercritical case}
\label{sec:globalexistence}

In this section we prove Theorem~\ref{globalexistence} in the supercritical case \(p>1+ \frac{\sigma}{N}\) by closing a bootstrap argument based on the energy inequality \eqref{eq:closed_ineq}.

\begin{proof}[Proof of Theorem~\ref{globalexistence}]
Define the integrating factor \(\mu(s) = \exp\bigl( - C \int_{s_0}^s e^{-\lambda_0\tau}\,d\tau \bigr)\), where \(C>0\) denotes the constant in the linear term of \eqref{eq:closed_ineq}. Since \(\mu(s_0)=1\) and \(\lambda_0>0\), \(\mu(s)\) is uniformly positive and bounded:
\[
0<c_1 \le \mu(s) \le c_2 \le 1.
\]
Multiplying \eqref{eq:closed_ineq} by \(\mu(s)\) and integrating over \([s_0,s]\) gives
\begin{equation}\label{eq:muE}
\mu(s)E_5(s) \le C_1 E_5(s_0) + C_2 \int_{s_0}^s e^{-2\lambda_1\tau} E_5(\tau)^p \,d\tau + C_2 \int_{s_0}^s e^{-\lambda_1\tau} E_5(\tau)^{(p+1)/2} \,d\tau .
\end{equation}
Let \(M(s)=\sup_{\tau\in[s_0,s]}E_5(\tau)\). Since the integrands are monotone in \(M\), and using \(\mu(s) \ge c_1\), we obtain for all \(s_0 \le s < S_0(\varepsilon)\) (where \(S_0(\varepsilon)\) is the maximal existence time in terms of \(s\)):
\begin{equation}\label{eq:M_integral}
M(s) \le C_1 E_5(s_0) + C_3 \left( \int_{s_0}^s e^{-2\lambda_1\tau} \, d\tau \right) M(s)^p + C_3 \left( \int_{s_0}^s e^{-\lambda_1\tau} \, d\tau \right) M(s)^{(p+1)/2},
\end{equation}
where \(C_1, C_3 > 0\) are independent of \(\varepsilon\) and \(s\).

In the supercritical case, \(\lambda_1 > 0\), so the time integrals in \eqref{eq:M_integral} are uniformly bounded as \(s \to \infty\) by 
\[
K := \max \left\{ \int_{s_0}^\infty e^{-2\lambda_1\tau}\,d\tau,\ \int_{s_0}^\infty e^{-\lambda_1\tau}\,d\tau \right\} < \infty.
\]
Hence \eqref{eq:M_integral} simplifies to the algebraic inequality:
\begin{equation}\label{eq:algebraic}
M(s) \le C_1 E_5(s_0) + C_4 \bigl( M(s)^p + M(s)^{(p+1)/2} \bigr), \qquad C_4 = C_3 K.
\end{equation}
By the local existence theory and scaling transformation, \(E_5(s_0) \le C_0 \varepsilon^2 I_0\) with \(I_0 = \|(u_0, u_1)\|_{H^{\sigma/2, m} \times H^{0, m}}^2\). Absorbing constants, \eqref{eq:algebraic} becomes
\[
M(s) \le C_1 \varepsilon^2 I_0 + C_4 \bigl( M(s)^p + M(s)^{(p+1)/2} \bigr).
\]

Set
\[
H(x) = x - C_3 \bigl( x^p + x^{(p+1)/2} \bigr).
\]
Since \(p > 1\), \(H(0) = 0\) and \(H'(0) = 1 > 0\). Hence there exists \(x^* > 0\) such that \(H\) is strictly increasing on \([0, x^*]\) and
\[
H(x) \ge \frac{x}{2}, \qquad 0 \le x \le x^*.
\]
Choose \(\varepsilon_0 > 0\) so small that for all \(0 < \varepsilon \le \varepsilon_0\),
\[
C_1 \varepsilon^2 I_0 < H(x^*), \qquad 2C_1 \varepsilon^2 I_0 < x^*.
\]
By continuity and the strict monotonicity of \(H\) on \([0, x^*]\), the equation \(H(x) = C_1 \varepsilon^2 I_0\) has a unique root \(x_1 \in (0, 2C_1 \varepsilon^2 I_0]\). Since \(M(s_0) \le C_1 \varepsilon^2 I_0 < x_1\) and \(M(s)\) is continuous, the inequality \(M(s) \le x_1\) persists for all \(s \ge s_0\). 

Indeed, if \(M(s)\) ever exceeded \(x_1\), there would be a first time \(s_1>s_0\) such that \(M(s_1)=x_1\) and \(H(M(s_1))=C_1 \varepsilon^2 I_0\). But \(M(s_1)=x_1<x^*\) and \(H\) is strictly increasing on \([0,x^*]\), so \(H(M(s_1))>H(x_1)\) leads to a contradiction.

Returning to the original variables, the uniform bound on \(E_5\) prevents blow-up of the \(H^{\sigma/2,m}\times H^{0,m}\) norm, and therefore the solution extends globally in time. In particular, the lifespan \(T_\varepsilon = +\infty\) for all \(\varepsilon \in (0,\varepsilon_0]\).
\end{proof}

\section{Lower bound of the lifespan in the subcritical and critical cases}

In this section we derive the lower bounds of the lifespan in the subcritical and critical cases by analyzing the integral inequality \eqref{eq:M_integral} established in Proof of Theorem~\ref{globalexistence}.

\begin{proof}[Proof of Theorem~\ref{blowupsub} and Theorem~\ref{blowupcritical}]

From the local initial bound \(E_5(s_0) \le C_0 \varepsilon^2 I_0\), recall from \eqref{eq:M_integral} that \(M(s) = \sup_{\tau \in [s_0, s]} E_5(\tau)\) satisfies
\begin{equation}\label{eq:M_lifespan}
M(s) \le C_1 \varepsilon^2 I_0 + C_3 \left(\int_{s_0}^s e^{-2\lambda_1 \tau} d\tau\right) M(s)^p + C_3 \left(\int_{s_0}^s e^{-\lambda_1 \tau} d\tau\right) M(s)^{\frac{p+1}{2}}
\end{equation}
for all \(s_0 \le s < S_0(\varepsilon)\).

To estimate the growth of \(M(s)\), let \(S_1 = S_1(\varepsilon) \ge s_0\) be the first time such that
\[
M(S_1) = 2 C_1 \varepsilon^2 I_0.
\]
If \(S_1 = \infty\), the solution is global and the lower bound holds trivially. Assuming \(S_1 < \infty\), substituting \(s = S_1\) into \eqref{eq:M_lifespan} yields
\[
2 C_1 \varepsilon^2 I_0 \le C_1 \varepsilon^2 I_0 + C_3 \left(\int_{s_0}^{S_1} e^{-2\lambda_1 \tau} d\tau\right) (2 C_1 \varepsilon^2 I_0)^p + C_3 \left(\int_{s_0}^{S_1} e^{-\lambda_1 \tau} d\tau\right) (2 C_1 \varepsilon^2 I_0)^{\frac{p+1}{2}},
\]
which simplifies to
\begin{equation}\label{eq:at_S1}
C_1 \varepsilon^2 I_0 \le C_5 \varepsilon^{2p} \int_{s_0}^{S_1} e^{-2\lambda_1 \tau} d\tau + C_5 \varepsilon^{p+1} \int_{s_0}^{S_1} e^{-\lambda_1 \tau} d\tau.
\end{equation}

\noindent \textbf{Case I: Subcritical condition} (\(1 < p < 1 + \frac{\sigma}{N}\)). 
Here \(\lambda_1 = -\gamma < 0\) with \(\gamma = \frac{N}{\sigma}(1 + \frac{\sigma}{N} - p) > 0\). For small \(\varepsilon > 0\), \(S_1\) is sufficiently large, and the integrals in \eqref{eq:at_S1} satisfy
\[
\int_{s_0}^{S_1} e^{k\gamma \tau} d\tau \ge \frac{1}{2k\gamma} e^{k\gamma S_1}, \qquad k = 1, 2.
\]
Thus, \eqref{eq:at_S1} yields
\[
\varepsilon^2 \le C \left( \varepsilon^{2p} e^{2\gamma S_1} + \varepsilon^{p+1} e^{\gamma S_1} \right).
\]
Since \(\varepsilon > 0\) is small, to prevent the right-hand side from being strictly smaller than the left-hand side, we must have
\[
e^{S_1} \ge C \varepsilon^{-\frac{p-1}{\gamma}}.
\]
Since the solution exists on \([s_0, S_1]\), the maximal existence time \(T_\varepsilon\) satisfies \(s(T_\varepsilon) \ge S_1\). Recalling that \(e^{s(T_\varepsilon)} = B(T_\varepsilon) + 1\) and that the map \(s \mapsto B(T)\) is monotonically increasing, we obtain the lower bound for the lifespan:
\begin{equation}
    B(T_\varepsilon) \ge e^{S_1} - 1 \ge C \varepsilon^{-\frac{\sigma(p-1)}{N + \sigma - Np}}.
\end{equation}

\noindent \textbf{Case II: Critical condition} (\(p = 1 + \frac{\sigma}{N}\)). 
Here, \(\lambda_1 = 0\). The integrals in \eqref{eq:at_S1} grow linearly: 
\[\int_{s_0}^{S_1} d\tau = S_1 - s_0 \sim S_1. \]
Thus, \eqref{eq:at_S1} reduces to
\[
C_0' \varepsilon^2 I_0 \le C_2' (S_1 - s_0) \varepsilon^{2p} + C_2' (S_1 - s_0) \varepsilon^{p+1}.
\]
Since \(p > 1\), we have \(2p > p+1\), which implies \(\varepsilon^{2p} \ll \varepsilon^{p+1}\) for small \(\varepsilon\). Therefore, the \(\varepsilon^{p+1}\) term dominates the right-hand side, leading to
\[
\varepsilon^2 \le C (S_1 - s_0) \varepsilon^{p+1} \implies S_1 - s_0 \ge C \varepsilon^{-(p-1)} = C \varepsilon^{-\frac{\sigma}{N}}.
\]
Again, since the solution exists on \([s_0, S_1]\), we have \(s(T_\varepsilon) \ge S_1\). Reverting to \(T_\varepsilon\) using \(s(T_\varepsilon) = \log(B(T_\varepsilon) + 1)\), this yields the exponential lower bound:
\begin{equation}
    B(T_\varepsilon) \ge \exp(S_1) - 1 \ge \exp\left( C \varepsilon^{-\frac{\sigma}{N}} \right) - 1.
\end{equation}
This completes the proof.
\end{proof}

\section*{Acknowledgements}
The first author is supported by Science Foundation of Zhejiang Sci-Tech University (ZSTU) under Grant No.13062120-Y. The second author is supported by JSPS, the Grant-in-Aid for Scientific Research (C) (No. 23K03174) and the Grant-in-Aid for Transformative Research Areas (B) (No. 25K24910).

\end{document}